\documentclass[12pt,dvips]{article}
\usepackage{amstex}
\usepackage{epsf}

\def\CC{{\rm\kern.24em\vrule width.02em height1.4ex depth-.05ex\kern-.26em C}}
\def\QQ{{\rm\kern.24em\vrule width.02em height1.4ex depth-.05ex\kern-.26em Q}}
\def\PP{{\rm\kern.24em\vrule width.02em height1.4ex depth-.05ex\kern-.26em P}}
\def\Rr{{\rm I\kern-.2em R}}
\def\ZZ{{\rm\kern.26em\vrule width.02em height0.5ex depth0ex\kern.04em\vrule width.02em height1.47ex depth-1ex\kern-.34em Z}}
\def\BB{{\rm\kern.24em\vrule width.02em height1.4ex depth-.05ex\kern-.26em B}}
\def\RR{{\hspace{.065in}\rm{\vrule width.02em height1.55ex depth-.07ex\kern-.3165em R}}}
\def\Ibb#1{{\rm I\kern-.23em#1}}
\def\NN{\Ibb N}

\newcommand {\interior} {\rm int \,}

\newcommand {\be} {\begin{equation}}
\newcommand {\ee} {\end{equation}}
\newcommand {\ba} {\begin{align}}
\newcommand {\ea} {\end{align}}
\def\bas{\begin{align*}}
\def\eas{\end{align*}}

\newcommand {\ra} {\rightarrow}

\newcommand {\BOX} { \hfill \rule{2mm}{2mm}}

\newtheorem {lemma} {LEMMA} [section]
\newtheorem {theorem}[lemma]{THEOREM} 
\newtheorem {prop}[lemma]{PROPOSITION}
\newtheorem {cor}[lemma]{COROLLARY}
\newtheorem {definition}[lemma]{Definition}
\newtheorem {exercise}[lemma]{Exercise}
\newtheorem {example}[lemma]{Example}
\newtheorem {question}[lemma]{Question}
\newtheorem {problem}[lemma]{Problem}

\numberwithin{equation}{section}

\begin{document}

\title{\bf Complex dynamics in several variables}
\author{Lectures by John Smillie\\ Notes by Gregery T. Buzzard}
\date{}
\maketitle

\tableofcontents
\newpage

\begin{abstract}
These notes are the outgrowth of a series of lectures given at MSRI in
January 1995 at the beginning of the special semester in complex
dynamics and hyperbolic geometry.  In these notes, the primary aim is
to motivate the study of complex dynamics in two variables and to
introduce the major ideas in the field.  Hence the treatment of the
subject is mostly expository. 
\end{abstract}

\renewcommand{\thefootnote}{}
\footnote{Research at MSRI is 
   supported in part by NSF grant DMS-9022140}

\small
\normalsize

\section{Introduction}

These notes are the outgrowth of a series of lectures given at MSRI in
January 1995 at the beginning of the special semester in complex
dynamics and hyperbolic geometry.  The goal of these lectures was to
provide an introduction to the relevant ideas and problems of complex
dynamics in several variables and to provide a foundation for further
study and research.  There were parallel sessions in complex dynamics
in one variable, given by John Hubbard, and in hyperbolic geometry,
given by James Cannon, and notes for those series should also be
available.

In these notes, the primary aim is to motivate the study of complex
dynamics in two variables and to introduce the major ideas in the
field.  Hence the treatment of the subject is mostly expository.

\section{Motivation}

The study of complex dynamics in several variables can be motivated in
at least two natural ways.  The first is by analogy with the fruitful
study of complex dynamics in one variable.  However, since this latter
subject is covered in detail in the parallel notes of John Hubbard, we
focus here on the motivation coming from the study of real dynamics.

A classical problem in the study of real dynamics is the $n$-body
problem, which was studied by Poincar\'e.  For instance, we can think
of $n$ planets moving in space.  For each planet, there are 3
coordinates giving the position and 3 coordinates giving the velocity,
so that the state of the system is determined by a total of $6n$ real
variables.  The evolution of the system is governed by Newton's laws,
which can be expressed as a first order ordinary differential
equation.  In fact, the state of the system at any time determines the
entire future and past evolution of the system.  

\bigskip

To make this a bit more precise, set $k= 6n$.  Then the behavior of
the $n$ planets is modeled by a differential equation
$$ (\stackrel{.}{x}_1, \ldots, \stackrel{.}{x}_k) = F(x_1, \ldots,
x_k)$$
for some $F: \RR^k \ra \RR^k$.  Here $\stackrel{.}{x}$ denotes the
derivative of $x$ with respect to $t$.  

From the elementary theory of ODE's, we
know that this system has a unique solution $t \mapsto \phi_t(x_1,
\ldots, x_k)$ satisfying $\stackrel{.}{\phi} = F(\phi)$ and
$\phi_0(x_1, \ldots, x_k) = (x_1, \ldots, x_k)$.  

For purposes of studying dynamics, we would like to be able to say
something about the evolution of this system over time, given some
initial data.  That is, given $p \in \RR^k$, we would like to be able
to say something about $\phi_t(p)$ as $t$ varies.  For instance, a
typical question might be the following.

\begin{question} 
Given $p = (x_1, \ldots, x_k)$, is $\{\phi_t(p):
t \geq 0\}$ bounded?
\end{question}

Unfortunately, the usual answer to such a question is ``I don't
know.''  Nevertheless, it is possible to say something useful about
related questions, at least in some settings.  For instance, one
related problem is the following.

\begin{problem}  
Describe the set
$$ K^+ := \{p \in \RR^k : \{ \phi_t(p) : t \geq 0\} \; \; {\rm is \;
\; bounded}\}. $$
\end{problem}

\medskip
Although this question is less precise and gives less specific
information than the original, it can still tell us quite a bit
about the behavior of the system.

\section{Iteration of maps}

In the above discussion, we have been taking the approach of fixing a
point $p \in \RR^k$ and following the evolution of the system over
time starting from this point.  An alternative approach is to think of
all possible starting points evolving simultaneously, then taking a
snapshot of the result at some particular instant in time.  

To make this more precise, assume that the solution
$\phi_t(p)$ exists for all time $t$ and all $p \in \RR^k$.  In this
case, for fixed $t$, the map $\phi_t : \RR^k \ra \RR^k$ is a
diffeomorphism of $\RR^k$ and satisfies the group property
$$ \phi_{s+t} = \phi_s \circ \phi_t $$
for any $s$ and $t$.

In order to make our study more tractable, we make two
simplifications.

\paragraph{Simplification 1:} Choose some number $\alpha > 0$, which
we will call the {\em 
period}, and define $f = \phi_\alpha$.  Then $f$ is a diffeomorphism
of $\RR^k$, and given $p \in \RR^k$, the group property of $\phi$
implies that
$$ \phi_{n\alpha}(p) = \phi_\alpha \circ \cdots \circ \phi_\alpha(p) =
f^n(p).$$ 
That is, studying the behavior of $f$ under iteration is equivalent to
studying the behavior of $\phi$ at regularly spaced time intervals.

\paragraph{Simplification 2:} Set $k = 2$.  Although this
simplification means that we can 
no longer directly relate our model to the original physical problem,
the ideas and techniques involved in studying such a simpler model are
still rich enough to shed some light on the more realistic cases.
In fact, there are interesting questions in celestial mechanics which
reduce to questions about two dimensional diffeomorphisms, but here
we are focusing on the mathematical model rather than on the physical
system.

\bigskip
Finally, we make the following definition.

\begin{definition}
Given $p \in \RR^2$, set
\begin{align*}
{\cal O}^+(p) & :=  \{f^n(p) : n \geq 0\}, \\
{\cal O}^-(p) & :=  \{f^n(p) : n \leq 0\}, \\
{\cal O}(p) & :=  \{f^n(p) : n \in \ZZ\}.
\end{align*}
\end{definition}

With these simplifications and this definition, we can further
reformulate the question from the previous section as follows.

\begin{problem}  \label{prob:K+-}
Given a diffeomorphism $f: \RR^2 \ra \RR^2$, describe the sets 
\begin{align*}
K^+ & :=  \{ p \in \RR^2 : {\cal O}^+(p) \; \; {\rm is \;\; bounded}\}, \\
K^- & := \{ p \in \RR^2 : {\cal O}^-(p) \; \; {\rm is \;\; bounded}\}, \\
K & :=  \{ p \in \RR^2 : {\cal O}(p) \; \; {\rm is \;\; bounded}\},
\end{align*}
\end{problem}

For future reference, note that $K = K^+ \cap K^-$.

\section{Regular versus chaotic behavior}

For the moment, we will make no attempt to define rigorously what we
mean by regular or chaotic.  Intuitively, one should think of regular
behavior as being very predictable and as relatively insensitive
to small changes in the system or initial conditions.  On the other
hand, chaotic behavior is in some sense random and can change
drastically with only slight changes in the system or initial
conditions.  Here is a relevant quote from Poincar\'e on chaotic
behavior:

\begin{quotation}
{\em A very small cause, which escapes us, determines a considerable effect
which we cannot ignore, and we say that this effect is due to chance.}
\end{quotation}

We next give some examples to illustrate both kinds of behavior,
starting with regular behavior.  First we make some definitions.

\begin{definition} 
A point $p \in \RR^2$ is a {\bfseries periodic
point} if $f^n(p) = p$ for some $n \geq 1$.  The smallest such $n$ is
the {\bfseries period} of $p$.
\end{definition}

\begin{definition} 
A periodic point $p$ is {\bfseries hyperbolic} if
$(Df^n)(p)$ has no eigenvalues on the unit circle.
\end{definition}

\begin{definition} 
If $p$ is a hyperbolic periodic point and both
eigenvalues are inside the unit circle, then $p$ is called a {\bfseries
sink}.
\end{definition}

\begin{definition} 
If $p$ is a hyperbolic periodic point, then
the set
$$ W^s(p) = \{q \in \RR^2: d(f^nq, f^np) \ra 0 \;\; {\rm as} \;\; n
\ra \infty\}$$
is called the {\bfseries basin of attraction} of $p$ if $p$ is
a sink, and the {\bfseries stable manifold} of $p$ otherwise.  Here $d$ is a
distance function.  If $f$ is a diffeomorphism and $p$ is not a sink,
then the {\bfseries unstable manifold} of $p$ is defined by replacing $f$ by
$f^{-1}$ in the above definition.
\end{definition}

\paragraph{Fact:} When $p$ is a sink, $W^s(p)$ is an open set
containing $p$.  

\bigskip
A sink gives a prime example of regular behavior.  Starting with any
point $q$ in the basin of attraction of a sink $p$, the forward orbit
of $q$ is asymptotic to the (periodic) orbit of $p$.  Since the basin is
open, this will also be true for any point $q'$ near enough to $q$.
Hence we see the characteristics of predictability and stability
mentioned in relation to regular behavior.

For an example of chaotic behavior, we turn to a differential equation
studied by Cartwright and Littlewood in 1940.  This is the system
$$ \stackrel{..}{y} - k (1-y^2)\stackrel{.}{y} + y = b \cos(t). $$
Introducing the variable $x = \stackrel{.}{y}$, we can write this as a
first-order system
\begin{align*}
\stackrel{.}{y} & = x \\
\stackrel{.}{x} & = g(x,y,t),
\end{align*}
where $g$ is a function satisfying $g(x,y,t+2 \pi) = g(x,y,t)$.  This
system has a solution $\phi_t$ as before with $\phi_t:\RR^2 \ra \RR^2$
a diffeomorphism.  Although the full group property does not hold for
$\phi$ since $g$ depends on $t$, we still have $\phi_{s+t} = \phi_s
\circ \phi_t$ whenever $s = 2 \pi n$ and $t = 2 \pi m$ for integers
$n$ and $m$.  Hence we can again study the behavior of this system by
studying the iterates of the diffeomorphism $f=\phi_{2 \pi}$.  

Rather than study this system itself, we follow the historical
development of the subject and turn to a more easily understood
example of chaotic behavior which was motivated by this system of
Cartwright and Littlewood: the Smale horseshoe.

\section{The horseshoe map and symbolic dynamics}  \label{sec:horseshoe}

The horseshoe map was first conceived by Smale as a way of
capturing many of the features of the Cartwright-Littlewood map in a
system which is easily understood.  

\begin{figure}
\epsfysize=2in \centering \leavevmode \epsfbox{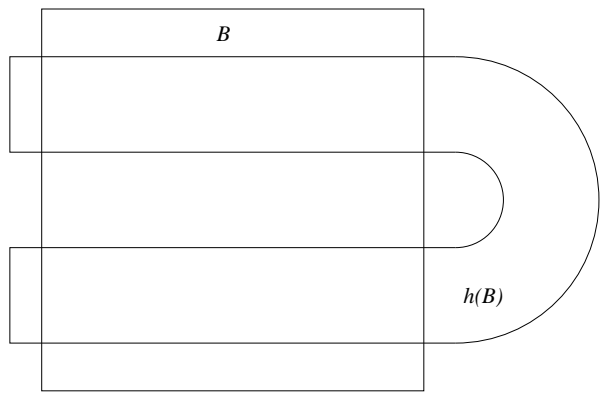}
\caption{The image of $B$ under the horseshoe map $h$.} \label{fig:h}
\end{figure}

For our purposes, the horseshoe map, $h$, is defined first on a square $B$ in
the plane with sides parallel to the axes.  First we apply a linear
map which stretches the square in the $x$-direction and contracts it
in the $y$-direction.  Then we take the right edge of the resulting
rectangle and lift it up and around to form a horseshoe shape.  The
map $h$ is then defined on $B$ by placing this horseshoe over the
original square $B$ so that $B \cap h(B)$ consists of 2 horizontal
strips in $B$.  See figure~\ref{fig:h}.

\begin{figure}
\begin{center}
\epsfysize=2.5in \leavevmode \epsfbox{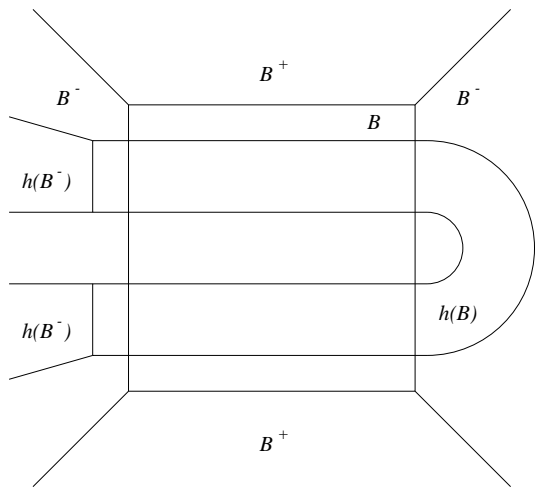}
\end{center}
\caption{The sets $B$, $B^+$, and $B^-$.} \label{fig:B}
\end{figure}

We can extend $h$ to a diffeomorphism of $\RR^2$ in many ways.  We do
that here as follows.  First partition $\RR^2 - B$ into 4 regions by
using the lines $y=x$ and $y=-x$ as boundaries.  Denote the union of
the two regions above and below $B$ by $B^+$ and the union of the
two regions to the left and right of $B$ by $B^-$, as in
figure~\ref{fig:B}.  Then we can 
extend $h$ to a diffeomorphism of $\RR^2$ in such a way that $h(B^-)
\subseteq B^-$.  In this situation, points in $B^+$ can be mapped to
any of the 3 regions $B^+$, $B$, or $B^-$, points in $B$ can be mapped
to either $B$ or $B^-$, and points in $B^-$ must be mapped to $B^-$.
Further, we require that points in $B^-$ go to $\infty$ under
iteration, and we require analogous conditions on $f^{-1}$.
Note in particular that points which leave $B$ do not return and that
$K \subseteq B$.

It is not hard to see that in this case, we have
$$ K^- \cap B = B \cap h B \cap h^2 B \cap \cdots.$$
In fact, if we look at the image of the two strips $B \cap h B$ and
intersect with $B$, then the resulting set consists of 4 strips; each of
the original two strips is subdivided into two smaller strips.
Continuing this process, we see that $K^- \cap B$ is simply the set
product of an interval and a Cantor set.  

In fact, a simple argument shows that $h$ has a fixed point $p$ in the
lower left corner of $B$, and that the unstable manifold of $p$ is dense
in the set $K^- \cap B$ and the stable manfold of $p$ is dense in $K^+
\cap B$.  The complicated structure of the stable and unstable
manifolds plays an important role in the behavior of the horseshoe map.

\bigskip
We can describe the chaotic behavior of the horseshoe using {\em symbolic
dynamics}.  The idea of this procedure is to translate from the dynamics of
$h$ restricted to $K$ into the dynamics of a shift map on bi-infinite
sequences of symbols.  

To do this, first label the 2 components of $B \cap hB$ with $H_0$ and
$H_1$.  Then to a point $p \in K$, we associate a bi-infinite sequence
of $0$'s and $1$'s using the map
$$ \psi:p \mapsto s=(\ldots, s_1, s_0, s_{-1}, \ldots), $$
where
$$ s_j = \left\{ 
\begin{array}{cl}
0 & {\rm if} \; \; h^j(p) \in H_0 \\
1 & {\rm if} \; \; h^j(p) \in H_1. 
\end{array}
\right.
$$

We can put a metric on the space of bi-infinite sequences of $0$'s and
$1$'s by 
$$ d(s,s') = \sum_{j=-\infty}^{\infty} |s_j-s_j'| \;2^{-|j|}.  $$
It is not hard to show that the metric space thus obtained
is compact and that the map $\psi$ given above 
produces a homeomorphism between $K$ and this space of sequences.

Moreover, the definition of $\psi(p)$ implies that if $\sigma$ is the
left-shift map defined on bi-infinite sequences, then 
$\psi(h(p)) = \sigma(\psi(p))$.  

Here are a couple of simple exercises which illustrate the power of
using symbolic dynamics.

\begin{exercise}  
Show that periodic points are dense in $K$.
Hint: Periodic points correspond to periodic sequences.
\end{exercise}

\begin{exercise}  \label{exs:per-pt}
Show that there are periodic points of all periods.
\end{exercise}

\section{H\'enon maps}

The horseshoe map was one motivating example for what are known as Axiom A
diffeomorphisms.  These received a great deal of attention in the 60's
and 70's.  Much current work focuses either on how Axiom A fails, as
in the work of Newhouse, or on how some Axiom A
ideas can be applied in new settings, as in the work of Benedicks and
Carleson \cite{bc} or Benedicks and Young \cite{by}. For more
information and further 
references, the book by Ruelle \cite{ruelle} provides a fairly gentle
introduction, while the books by Palis and de Melo \cite{pd}, Shub
\cite{shub}, and Palis and Takens \cite{pt} are more advanced.  See
also the paper by Yoccoz in \cite{y}.

Although there are still some interesting open questions about Axiom A
diffeomorphisms and related subjects, much of the study in the post
Axiom A era has centered around the study of the so-called H\'enon
map.  This is actually a family of diffeomorphisms $f_{a,b}: \RR^2 \ra
\RR^2$ defined by
$$ f_{a,b}(x,y) = (-x^2 + a - b y, x) $$
for $b \neq 0$.  These maps arise from a simplification of a
simplification of a map describing turbulent fluid flow.  

\begin{figure}
\begin{center}
\epsfysize=2in \leavevmode \epsfbox{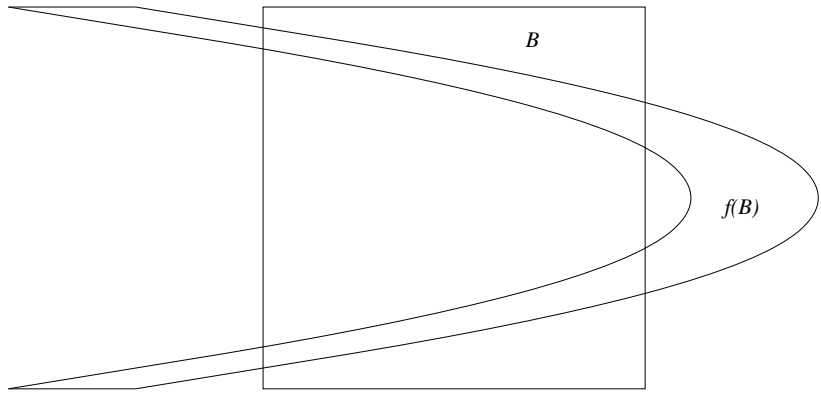}
\end{center}
\caption{A square $B$ and its image $f(B)$ for some parameter values
$a$ and $b$.}  \label{fig:f}
\end{figure}

\bigskip
We can get some idea of the behavior of the map $f_{a,b}$ and the ways
in which it relates to the horseshoe map by considering the image of a
large box $B$ under $f_{a,b}$.  For simplicity, we write $f$ for
$f_{a,b}$.  From figure~\ref{fig:f}, we see that for some
values of $a$ and $b$, the H\'enon map $f$ is quite reminiscent of the
horseshoe map $h$.

Since the map $f$ is polynomial in $x$ and $y$, we can also
think of $x$ and $y$ as being complex-valued.  In this case, $f:
\CC^2 \ra \CC^2$ is a holomorphic diffeomorphism of $\CC^2$.  This is
also in some sense a change in the map $f$, but all of the
dynamics of $f$ restricted to $\RR^2$ are contained in the
dynamics of the maps on $\CC^2$, so we can still learn about the
original map by studying it on this larger domain.

\begin{figure}[t]
\begin{minipage}[t]{.5\linewidth}
 \centering\epsfysize=2in \leavevmode \epsfbox{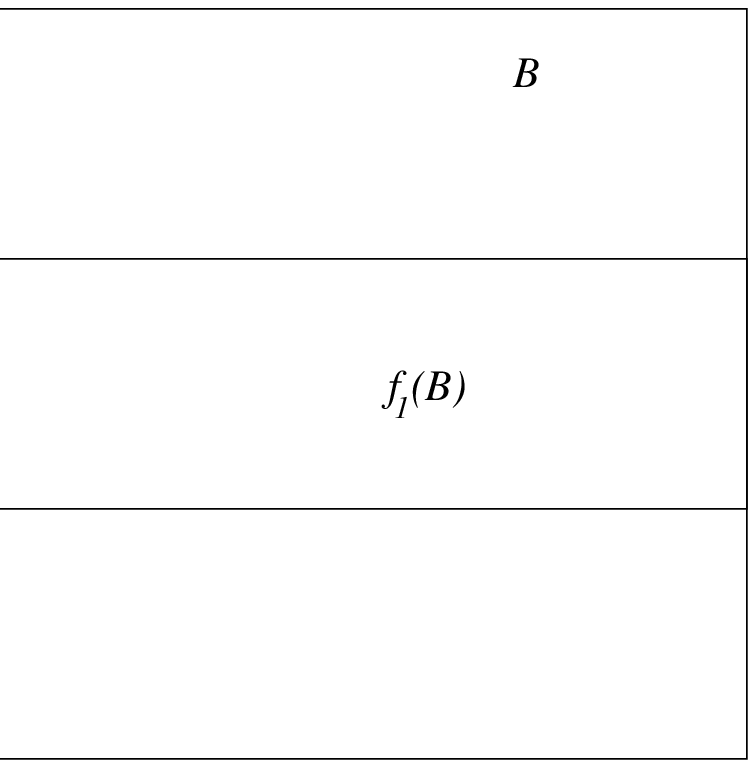}
 \caption{$f_1(B)$}   \label{fig:f1}
\end{minipage}
\begin{minipage}[t]{.49\linewidth}
 \centering\epsfysize=2in \leavevmode \epsfbox{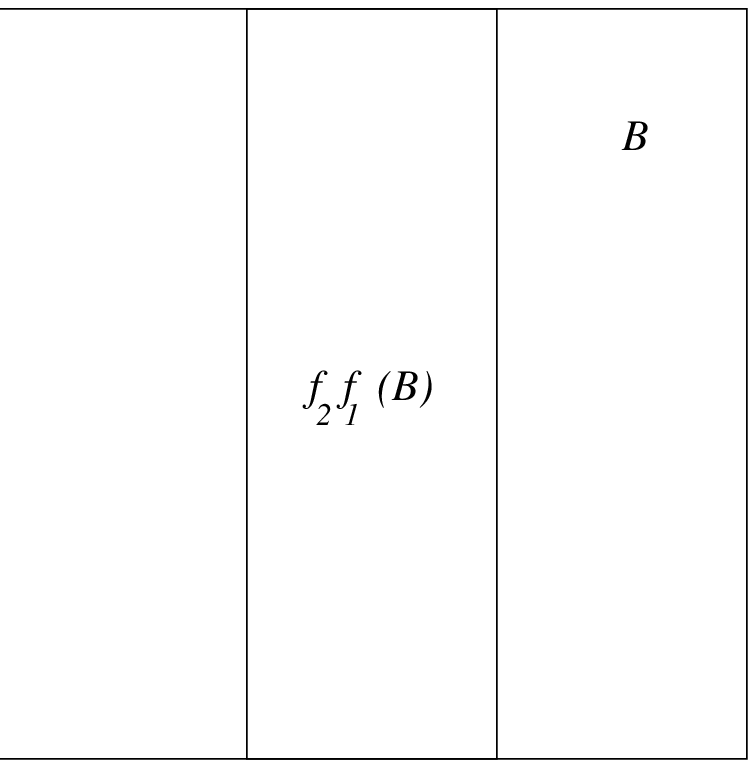}
 \caption{$f_2 f_1(B)$}  \label{fig:f2f1}
\end{minipage}
\end{figure}

\bigskip 
We next make a few observations about $f$.  First, note that $f$ is the
composition $f = f_3 \circ f_2 \circ f_1$ of the three maps
\begin{align*}
f_1(x,y) & = (x,by) \\
f_2(x,y) & = (-y,x) \\
f_3(x,y) & = (x + (-y^2 + a), y)
\end{align*}
For $0<b<1$, the maps $f_1$ and $f_2 f_1$ are depicted in
figures~\ref{fig:f1} and \ref{fig:f2f1}, while $f$ is depicted in
figure~\ref{fig:f} with some $a>0$. 

From the composition of these functions, we can easily see that $f$
has constant Jacobian determinant $\det(DF) = b$.  Moreover, when
$b=0$, $f$ reduces to a quadratic polynomial on $\CC$.

\bigskip
A simple argument shows that there is an $R = R(a,b)$ such that if we
define the three sets
\begin{align*}
B & = \{|x|<R, |y|<R\}, \\
B^+ & = \{|y|>R, |y|>|x|\}, \\
B^- & = \{|x|>R, |x|>|y|\},
\end{align*}
then we have the same dynamical relations as for the corresponding sets for
the horseshoe map.  That is, points in $B^+$ can be mapped to $B^+$,
$B$, or $B^-$, points in $B$ can be mapped to $B$ or $B^-$, and points
in $B^-$ must be mapped to $B^-$.  

As we did for the horseshoe, we can define $K^+$ to be the set of
points with bounded forward orbit, $K^-$ to be the set of points with
bounded backward orbit, and $K$ to be the intersection of these two sets.  

When $a$ is large enough, so that the tip of $f_{a,b}(B)$ is
completely outside $B$, then $f_{a,b}$ ``is'' a horseshoe.  By this we
mean that $B \cap f_{a,b}$ has two components, and hence we can use
symbolic dynamics exactly as before to show that the dynamics of
$f$ restricted to $K$ are exactly the same as in the horseshoe case.
The general proof of this result is due to Hubbard and Oberste-Vorth
\cite{ho}.  See also \cite{dn}.

\begin{example}  
To compare the dynamics of $f$ in the real and
complex cases, consider $f_{a,b}$ with $a$ and $b$ real.  As an ad hoc
definition, let $K_{\bf R}$ be the set of $p \in \RR^2$ with bounded
forward and backward orbits, and let $K_{\bf C}$ be the set of $p \in
\CC^2$  with bounded forward and backward orbits.  Then another result
of \cite{ho} is that $K_{\bf C} = K_{\bf R}$.
\end{example}

Thus we already have a mental picture of $K$ for these parameter
values.  We can also get a picture of $K^+$ and $K^-$ in the complex
case, since we can extend the analogy between $f$ and the 
horseshoe map by replacing the square $B$ by a bidisk $B = D(R) \times
D(R)$ contained in $\CC^2$, where $D(R)$ is the disk of radius $R$
centered at $0$ in $\CC$.  In the definitions of $B^+$ and $B^-$, we
can interpret $x$ and $y$ as complex-valued, in which case the
definitions of these sets still make sense.  Moreover, the same
mapping relations hold among $B^+$, $B^-$, and $B$ as before.  In this
case, $B \cap K^+$ is topologically equivalent to the 
set product of a Cantor set and a disk, $B \cap K^-$ is equivalent to
the product of a disk and a Cantor set, and $B \cap K$ is equivalent
to the product of two Cantor sets.

\paragraph{Thesis:}  A surprising number of properties of the
horseshoe (when properly interpreted) hold for general complex H\'enon
diffeomorphisms. 

\bigskip
The ``surprising'' part of the above thesis is that the horseshoe map was
designed to be simple and easily understood, yet it sheds much light
on the less immediately accessible H\'enon maps.

\section{Properties of horseshoe and H\'enon maps}

We again consider some properties of the horseshoe map in terms of
its periodic points.  There is another relevant
quote from Poincar\'e about periodic points.

\begin{quotation}
{\em What renders these periodic points so precious to us is that they
are, so to speak, the only breach through which we might try to
penetrate into a stronghold hitherto reputed unassailable.}
\end{quotation}

As an initial
observation, recall that from symbolic dynamics, we know that the
periodic points are dense in $K$.  In fact, it is not hard to show
that these periodic points are all {\em\bfseries saddle points}; that
is, if $p$ has 
period $n$, then $(D h^n)(p)$ has one eigenvalue larger than 1 in
modulus, and one smaller.  For such a periodic point $p$, recall the
definitions of the stable and unstable manifolds
\begin{align*}
W^s(p) & = \{q: d(f^np, f^nq) \ra 0 \;\; {\rm as} \;\; n \ra \infty\} \\
W^u(p) & = \{q: d(f^np, f^nq) \ra 0 \;\; {\rm as} \;\; n \ra -\infty\}.
\end{align*}

\begin{exercise}  \label{exs:density}
For any periodic saddle point of the horseshoe
map $h$, $W^s(p)$ is dense in $K^+$ and $W^u(p)$ is dense in $K^-$.
\end{exercise}

As a consequence of this exercise, suppose $p \in K^+$, and let $n \in
\NN$ and $\epsilon > 0$.  By exercise~\ref{exs:per-pt}, there is a periodic
point $q$ with period $n$, and by this last exercise, the stable
manifold for $q$ comes arbitrarily close to $p$.  In particular, we
can find $p' \in W^s(q)$ with $d(p, p') < \epsilon$.  Hence in any
neighborhood of $p$, there are points which are asymptotic to a
periodic point of any given period.  We can contrast this with a point $p$
in the basin of attraction for a sink.  In this case, for a small
enough neighborhood of $p$, every point will be asymptotic to the same
periodic point.  

This example illustrates the striking difference between regular and
chaotic behavior.  In the case of a sink, the dynamics of the map are
relatively insensitive to the precise initial conditions, at least
within the basin of attraction.  But in the horseshoe case, the
dynamics can change dramatically with an arbitrarily small change in
the initial condition.  In a sense, chaotic behavior occurs throughout
$K^+$.

A second basic example of Axiom A behavior is the solenoid.  Take a
solid torus in $\RR^3$ and map it inside itself so that it wraps
around twice.  The image of this new set then wraps around 4 times.
The solenoid is the set which is the intersection of all the
forward images of this map.  Moreover, the map extends to a
diffeomorphism of $\RR^3$ and displays chaotic behavior on the
solenoid, which is the attractor for the diffeomorphism.

\begin{example}  Consider $f_{a,b}: \CC^2 \ra \CC^2$ when $a$ and
$b$ are small.  Then Hubbard and Oberste-Vorth \cite{ho} show that
$f_{a,b}$ has a fixed sink as well as a solenoid, so that it
displays both regular and chaotic behavior.  
\end{example}

Note that if $q$ is a sink, then $W^s(q) \subseteq K^+$ is open, and hence
$W^s(q) \subseteq \interior K^+$.  On the interior of $K^+$, there is
no chaos.  To see this, suppose $p \in \interior K^+$, and choose
$\epsilon > 0$ such that $\overline{\BB_\epsilon(p)} \subseteq K^+$.
A simple argument using the form of $f$ and the definitions of $B$,
$B^+$, and $B^-$ shows that any point in $K^+$ must eventually be
mapped into $B$.  Hence by compactness, there is an $n$ sufficiently
large that $f^n(\overline{\BB_\epsilon(p)}) \subseteq B$.  Since $B$
is bounded, we see by Cauchy's integral formula that the norm of the
derivatives of $f^n$ are uniformly bounded on $\BB^\epsilon(p)$
independently of $n \geq 0$.  This is incompatible with chaotic
behavior.  For more information and further references, see \cite{bs3}.

\bigskip
As a first attempt at studying sets where chaotic behavior can occur,
we make the following definitions.

\begin{definition} 
For a complex H\'enon map $f$, and with $K^+$ and $K^-$ defined as in
problem~\ref{prob:K+-}, let $J^+ := \partial K^+$ and $J^- := \partial K^-$.
\end{definition}

The following theorem gives an analog of exercise~\ref{exs:density} in
the case of a general complex H\'enon mapping, and is contained in
\cite{bs1}. 

\begin{theorem}
If $p$ is a periodic saddle point of the H\'enon map $f$, then $W^s(p)$
is dense in $J^+$, and $W^u(p)$ is dense in $J^-$.
\end{theorem}

It can be shown \cite{bls} that a H\'enon map $f$ has saddle periodic
points of all but finitely many periods, so just as in the argument
after exercise~\ref{exs:density}, we see that chaotic behavior occurs
throughout $J^+$, and a similar argument applies to $J^-$ under
backward iteration.

\section{Dynamically defined measures}  \label{sec:m-}

Before talking about potential theory proper, we first discuss some
measures associated with the horseshoe map $h$.  With notation as in
section~\ref{sec:horseshoe}, we define the level-$n$ set of $h$ to be
the set $h^{-n} B \cap h^n B$.  Since the forward images of $B$ are
horizontal strips and the backward images of $B$ are vertical strips,
we see that the level-$n$ set consists of $2^{2n}$ disjoint boxes.

\paragraph{Assertion:} For $j$ sufficiently large, the number of fixed
points of $h^j$ in a component of the level-$n$ set of $h$ is
independent of the component chosen.

\bigskip
In fact, there is a unique probability measure $m$ on $K$ which
assigns equal weight to each level-$n$ square, and the above assertion
can be rephased in terms of this measure.  Let $P_k$ denote the
set of $p \in \CC^2$ such that $h^k(p) = p$.  Then it follows from the
above assertion that 
\begin{equation}   \label{eqn:m} 
\frac{1}{2^k} \sum_{p \in P_k} \delta_p \ra m 
\end{equation}
in the topology of weak convergence.

\begin{figure}
\begin{center}
\epsfysize=2in \leavevmode \epsfbox{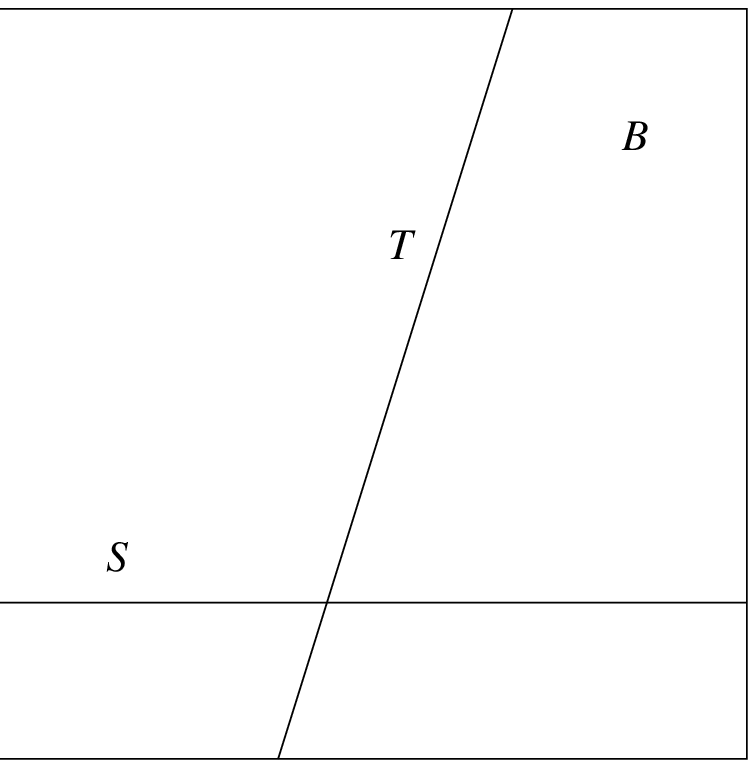}
\end{center}
\caption{$T$ and $S$}   \label{fig:ST1}
\end{figure}

\bigskip
We can use a similar technique to study the distribution of unstable
manifolds.  Again we consider the horseshoe map $h$, and we suppose
that $p_0$ is a fixed saddle point of $h$ and that $S$ is the component of
$W^u(p_0) \cap B$ containing $p_0$.  In this case, $S$ is simply a
horizontal line segment through $p_0$.  Next, let $T$ be a line
segment from the top to the bottom of $B$ so that $T$ is transverse to
every horizontal line.  See figure~\ref{fig:ST1}.

We can define a measure on $T$ using an averaging process as before.
This time we average over points in $h^n(S) \cap T$ to obtain a
measure $m_T^-$.  Thus we have
\begin{equation}   \label{eqn:m_T}
\frac{1}{2^n} \sum_{p \in h^n(S) \cap T} \delta_p \ra m_T^-, 
\end{equation}
where again the convergence is in the weak sense.  This gives a
measure on $T$ which assigns equal weight to each level-$n$ segment;
i.e., to each component of $h^n(B) \cap T$.

Note that if $T'$ is another segment like $T$, then the unstable
manifolds of $h$ give a way to transfer the above definition to a
measure on $T'$.  That is, given a point $p \in h^n(S) \cap T$, 
we can project along the component of $h^n(S) \cap B$ containing $p$
to obtain a point $p' \in T'$.  Using this map we obtain a measure
$\phi(m_T^-)$ on $T'$.  It is straightforward to show that this is the
same measure as $m_{T'}^-$ obtained by using $T'$ in place of $T$ in
(\ref{eqn:m_T}).  

Given this equivalence among these measures, we can define $m^-$ to be
this collection of equivalent measures, so that $m^-$ is defined on
any $T$.  Using an analogous construction with stable manifolds, we
can likewise define a measure $m^+$ defined on ``horizontal''
segments.  

Finally, we can take the product of these two measures to
get a measure $m^- \times m^+$ defined on $B$.  Then one can show that
this product measure is the same as the measure $m$ defined in
(\ref{eqn:m}).  Hence there are at least two dynamically natural ways
to obtain this measure.

\section{Potential theory}

In the study of dynamics in one variable, there are many tools
available coming from classical complex analysis, potential theory,
and the theory of quasiconformal mappings.  In higher dimensions, not
all of these tools are available, but one tool which remains useful is
potential theory.  In this section we will provide some background for
the ways in which this theory can be used to study dynamics.

To provide some physical motivation for the study of potential theory,
consider two electrons moving in $\RR^d$, each with a charge of $-1$.
Then the repelling force between them is proportional to $1/r^{d-1}$.
If we fix one electron at the origin, then the total work in moving
the other electron from the point $z_0$ to the point $z_1$ is
independent of the path taken and is given by $P(z_1) - P(z_0)$,
where $P$ is a potential function which depends on the dimension:
\begin{equation}  \label{eqn:P}
\begin{array}{rcll}
P(z) & = & |z| &  {\rm if} \; \; d=1,\\
P(z) & = & \log|z| & {\rm if} \; \; d=2,\\
P(z) & = & -\frac{1}{\|z\|^{d-2}} & {\rm if} \;\;d \geq 3.
\end{array}
\end{equation}

From the behavior of $P$ at $0$ and $\infty$, we see that if $d \leq
2$, then the amount of work 
needed to bring a unit charge in from the point at infinity is infinite,
while this work is finite for $d \geq 3$.  On the other hand, if $d
\geq 2$, then the amount of work needed to bring two electrons together
is infinite, but for $d=1$ this work is finite.

\bigskip
We can think of a collection of electrons as a charge, and we can
represent charges by measures $\mu$ on $\RR^d$.   
Then for $S \subseteq \RR^d$, $\mu(S)$ is the amount of charge on $S$.

\begin{example}
A unit charge at the point $z_0$  corresponds to the Dirac delta mass
$\delta_{z_0}$. 
\end{example}

By using measures to represent charges, we can use convolution to
define potential functions for general charge distributions.  That is,
given a measure $\mu$ on $\RR^d$, we define
\begin{equation}  \label{eqn:pmu}
P_\mu(z) = \int_{\bf R^d} P(z-w) d\mu(w), 
\end{equation}
where $P$ is the appropriate potential function from (\ref{eqn:P}).
Note that this definition agrees with the previous definition of
potential functions in the case of point charges.  Note also that the
assignment $\mu \mapsto P_\mu$ is linear in $\mu$.

\bigskip
In order to be able to use potential functions to study dynamics, we
first need to understand a little more about their properties.  In
particular, we would like to know which functions can be the potential
function of a finite measure.

In the case $d=1$, the definition of $P_\mu$ and the triangle
inequality imply that potential functions are convex, hence also
continuous.  We also have that $P_\mu(x) = c|x| +
O(1)$, where $c = \mu(\RR)$.  In fact, any function, $f$, satisfying
these two conditions is a potential function of some measure.  Hence a
natural question is, how do we recover the measure from $f$?

In particular, given a convex function $f$ of one real variable, we
can consider the assignment 
\begin{equation}  \label{eqn:sd}
f \mapsto \frac{1}{2}\left( \frac{\partial ^2}{\partial x^2} f
\right) dx, 
\end{equation}
where the right hand side is interpreted in the sense of
distributions.  By convexity, this distribution is positive.  That is,
it assigns a positive number to positive test functions, and a
positive distribution is actually a positive measure.  Hence we have
an explicit correspondence between 
convex functions and positive measures, and with the additional
restriction on the growth of potential functions given in the previous
paragraph, we have an explicit correspondence between potential
functions and finite positive measures. (The $1/2$ in the above
formula occurs because we have normalized by dividing by the ``volume'' of
the unit sphere in $\RR$, i.e., the volume of the two points $1$ and $-1$.)

\bigskip
In the case $d=2$, the integral definition of $P_\mu$ implies that
potential functions satisfy the {\em\bfseries sub-average property}.
That is, given 
a potential function $f$, any $z_0$ in the plane, and a disk $D$
centered at $z_0$, $f(z_0)$ is bounded above by the average of $f$ on
$\partial D$.  That is, if $\sigma$ represents 1-dimensional Legesgue
measure normalized so that the unit circle has measure 1, and if $r$
is the radius of $D$, then
$$ f(z_0) \leq \frac{1}{r} \int_{\partial D} f(\zeta) d
\sigma(\zeta). $$

Moreover, (\ref{eqn:pmu}) implies that potential functions are
{\em\bfseries upper 
semicontinuous} (u.s.c.); a real-valued function $f$ is said to be
u.s.c. if all of its sub-level sets are open.  With these two concepts
we can make the following definition.

\begin{definition}
If $f$ is u.s.c. and satisfies the sub-average property, then $f$ is
called {\bfseries subharmonic}.
\end{definition}

Finally, if $f$ is subharmonic and satisfies $f(z) = c \log|z| + O(1)$
for some $c>0$, then $f$ is said to be a {\em\bfseries potential function}.
Just as before, a potential function has the form $P_\mu$ for some
measure $\mu$.

\bigskip
In fact, if $f$ is subharmonic and $C^2$, then the Laplacian of $f$ is
always positive.  This is an analog of the fact that the second
derivative of a convex function is positive.  If $f$ is subharmonic
but not $C^2$, then $\Delta f$ is a positive distribution, hence a
positive measure.  Thus the Laplacian
gives us a correspondence between potential functions and finite
measures much like that in (\ref{eqn:sd}):
$$ f \mapsto \frac{1}{2 \pi} (\Delta f) dx dy,  $$
where this is to be interpreted in the sense of distributions and
again we have normalized by dividing by the volume of the unit 
sphere.  

\begin{example}
Applying the above assignment to the function $\log |z|$ produces the delta
mass $\delta_0$ in the sense of distributions.
\end{example}

Suppose now that $K \subseteq \RR^2 = \CC$ is compact and connected
and put a charge on $K$ and allow it to distribute evenly throughout $K$.
Then the charge on $K$ is distributed according to a finite positive
measure $\mu$, and we would like to know what the equilibrium state is
for this system. Thus, we want to know what $P_\mu$ looks like.  

To use more standard notation we write $G = P_\mu$, and we assume that
$\mu$ has total charge (mass) equal to 1.  Then $G$ satisfies the
following properties.

\begin{enumerate}
\item $G$ is subharmonic.
\item $G$ is harmonic outside $K$.
\item $G = \log|z| + O(1)$.
\item $G$ is constant on $K$.
\end{enumerate}
If $G$ satisfies properties 1 through 3, and also property
\begin{itemize}
\item[$4'$.]  $G \equiv 0$ on $K$,
\end{itemize}
then we say that $G$ is a Green function for $K$.  If $G$ exists, then
it is unique, and in this case we can take the Laplacian of $G$ in the
sense of distributions.  Thus, we say that
$$ \mu_K := \frac{1}{2\pi} \Delta G \; dx dy  $$
is the equilibrium measure for $K$.

\begin{example}
Let $D$ be the unit disk.  Then the Green function for $D$ is
$$ G(z) = \log^+ |z|, $$
where $\log^+|z| := \max \{ \log|z|, 0\}$, and the equilibrium measure
is 
$$ \mu_D = \frac{1}{2\pi} (\Delta \log^+|z|) dx dy, $$
which is simply arc length measure on $\partial D$, normalized to have
mass 1.
\end{example}

\section{Potential theory in one variable dynamics}

In this section we discuss some of the ways in which potential theory
can be used to understand the dynamics of holomorphic maps of the
Riemann sphere.  This idea was first introduced by Brolin
\cite{brolin} and later developed by others in both one and several
variables.  

\bigskip
For this section, let $f$ be a monic polynomial in one variable of
degree $d \geq 2$, and let $K \subseteq \CC$ be the set of $z$ such
that the forward orbit of $z$ is bounded.  Then $K$ has a Green
function, and in fact, $G_K$ is given by the formula
$$ G_K(z) = \lim_{n \ra \infty} \frac{1}{d^n} \log^+|f^n(z)|.  $$

It is difficult to understand Brolin's paper without knowing this
formula.  However, it was in fact first written down by Sibony in his
UCLA lecture notes after Brolin's paper had already been written.

It is not hard to show that the limit in the definition of $G_K$
converges uniformly on compact 
sets, and since each of the functions on the right hand side is
subharmonic, the limit is also subharmonic.  Moreover, on a given
compact set outside of $K$, each of these functions is harmonic for
sufficiently large $n$, so that the limit is harmonic on the
complement of $K$.  Property 3 follows by noting that for $|z|$ large
we have $|z|^d/c \leq |f(z)| \leq c |z|^d$ for some $c > 1$, then taking
logarithms and dividing by $d$, then using an inductive argument to
bound $|\;\log^+|f^n(z)|/d^n - \log|z|\;|$ independently of $d$.
Finally, property $4'$ is immediate since $\log^+|z|$ is bounded for
$z \in K$.  In fact, $G_K$ has the additional property of being continuous.

Hence we see that $G_K$ really is the Green function for $K$, and we
can define the equilibrium measure 
$$ \mu := \mu_K = \frac{1}{2 \pi} (\Delta G_K) dx dy.  $$

\bigskip
The following theorem provides a beautiful relationship between the
measure $\mu$ and the dynamical properties of $f$.  It says that we
can recover $\mu$ by taking the average of the point masses at the
periodic points of period $n$ and passing to the limit or by taking
the average of the point masses at the preimages of any nonexceptional
point and passing to the limit.  A point $p$ is said to be nonexceptional
for a polynomial $f$ if the set $\{f^{-n}(p): n \geq 0\}$ contains at
least 3 points.  It is a theorem that there is at most 1 exceptional
point for any polynomial.  

\begin{theorem}
[Brolin, Tortrat]  Let $f$ be a monic polynomial of degree $d$, and
let $c \in \CC$ be a nonexceptional point.  Then 
$$ \mu = \lim_{n \ra \infty} \frac{1}{d^n} \sum_{z \in A_n} \delta_z,
$$
in the sense of weak convergence, where $A_n$ is either \\[5mm]
1. The set of $z$ satisfying $f^n(z) = c$, counted with
multiplicity, \\[3mm]
or \\[3mm]
2. The set of $z$ satisfying $f^n(z) = z$, counted with
multiplicity.
\end{theorem}

{\bf Proof:}  We prove only part 1 here.  Let 
$$ \mu_n = \frac{1}{d^n} \sum_{f^n(z) = c} \delta_z.  $$
Then we want to show that $\mu_n \ra \mu_K$.  Since the space of
measures with the topology of weak convergence is compact, it suffices
to show that if some subsequence of $\mu_n$ converges to a measure
$\mu^*$, then $\mu^* = \mu_K$.  By renaming, we may assume that
$\mu_n$ converges to $\mu^*$.  

We can show $\mu^* = \mu_K$ by showing the convergence of the
corresponding potential functions.  The potential function for $\mu_n$
is 
\begin{align*}
G_n(z) & = \frac{1}{d^n} \sum_{f^n(w) = c} \log |z-w| \\
& = \frac{1}{d^n} \log \bigg| \prod_{f^n(w) - c = 0} (z-w)\bigg| \\
& = \frac{1}{d^n} \log|f^n(z)-c|.
\end{align*}
Here the sum and products are taken over the indicated sets with
multiplicities, and the last equality follows from the fact that we
are simply multiplying all the monomials corresponding to roots of
$f^n(z)-c$.  

\bigskip
Let $G^*(z) := \lim_{n \ra \infty} G_n(z)$.  Then $G^*$ is the
potential function for $\mu^*$, and 
$$ G^*(z) = \lim_{n \ra \infty} \frac{1}{d^n} \log |f^n(z) - c|, $$
while
$$ G_K(z) = \lim_{n \ra \infty} \frac{1}{d^n} \log^+ |f^n(z)|, $$
and we need to show that $G^*(z) = G_K(z)$.  If $z \not \in K$, then
$f^n(z)$ tends to $\infty$ as $n$ increases, so that $G^*(z) = G_K(z)$
in this case.  Since $G^*$ is the potential function for $\mu^*$, it
is upper semi-continuous, so it follows that $G^*(z) \geq 0$ for $z
\in \partial K$.  On the 
other hand, since $G^* = G$ on the set where $G = \epsilon$, the
maximum principle for subharmonic functions implies that $G \leq
\epsilon$ on the region enclosed by this set.  Letting $\epsilon$ tend
to $0$ shows that $G^* \leq 0$ on $K$.  

Finally, using some knowledge of the possible types of components for
the interior of $K$, it follows that if $c$ is non-exceptional, then
the measure $\mu^*$ 
assigns no mass to the interior of $K$.  This implies that $G^*$ is
harmonic on $K$ since $\mu^*$ is the Laplacian of $G^*$.  Hence both
the maximum and minimum principles apply to $G^*$ on $K$, which
implies that $G^* \equiv 0$ on $K$.  

Thus $G^* \equiv G_K$ and hence $\mu^* \equiv \mu_K$ as desired. $\BOX$

\paragraph{Remark:}  This theorem provides an algorithm for drawing
a picture of the Julia set, $J$, for a polynomial $f$.  That is, start with
a nonexceptional point $c$, and compute points on the backward orbits
of $c$.  These points will accumulate on the Julia set for $f$, and by
discarding points in the first several backwards iterates of $c$, we
can obtain a reasonably good picture of the Julia set.  This algorithm
has the disadvantage that these backwards orbits tend to accumulate
most heavily on points in $J$ which are easily accessible from
infinity. That is, it favors points at which a random walk starting at
infinity is most likely to land and avoids points such as inward
pointing cusps. 

\begin{exercise}
Since $G$ is harmonic both in the complement of $K$ and in the
interior of $K$, we see 
that \;${\rm supp} \; \mu \subseteq J$, where $J = \partial K$ is the Julia
set.  Show that ${\rm supp} \; \mu = J$.  Hint:  Use the maximum
principle. 
\end{exercise}

Note that an immediate corollary of this exercise and Brolin's theorem
is that periodic points are dense in $J$.

\section{Potential theory and dynamics in two variables}

In Brolin's theorem, we took the average of point masses distributed
over the following two sets:
\begin{itemize}
\item[(a)] $\{z: f^n(z) = c\}$
\item[(b)] $\{z: f^n(z) = z\}$.
\end{itemize}
In the setting of polynomial automorphisms of $\CC^2$, there are two
natural questions motivated by these results.
\begin{itemize}
\item[(i)] What happens when we iterate 1-dimensional
submanifolds (forwards or backwards)?
\item[(ii)] Are periodic points described by some measure $\mu$?
\end{itemize}

In $\CC$, we can loosely describe the construction of the measure
$\mu$ as first counting the number of points in the set (a) or (b),
then using potential theory to describe the location of these points.

\bigskip
Before we consider such a procedure in the case of question (i) for
$\CC^2$, we first return to the horseshoe map and recall the measure
$m^-$ defined in section~\ref{sec:m-}.  Suppose that $B$ is
defined as in that section, that $p$ is a fixed point for the
horseshoe map $h$, that $S$ is the component of $W^u(p) \cap B$
containing $p$, and that $T$ is a line segment from the top to the
bottom of $B$ as before.  Then orient $T$ and $S$ so that these
orientations induce the standard orientation on $\RR^2$ at the point
of intersection of $T$ and $S$.  

\begin{figure}
\begin{center}
\epsfysize=2in \leavevmode \epsfbox{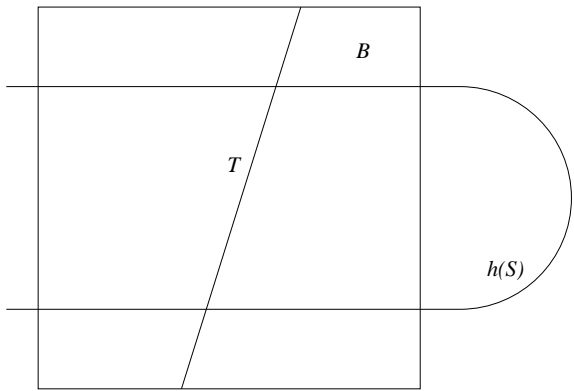}
\end{center}
\caption{$T$ and $h(S)$}   \label{fig:ST2}
\end{figure}

Now, apply $h$ to $S$.  Then $h(S)$ and $T$ will intersect in two
points, one of which is the original point of intersection, and one
of which is new.  See figure~\ref{fig:ST2}.  Because of the form of the
horseshoe map, the 
intersection of $h(S)$ and $T$ at the new point will not induce the
standard orientation on $\RR^2$ but rather the opposite orientation.
In general, we can apply $h^n$ to $S$, then assign $+1$ to each point
of intersection which induces the standard orientation, and $-1$ to
each point which induces the opposite orientation.  Unfortunately,
the sum of all such points of intersection for a given $n$ will always
be $0$, so this doesn't give us a way to count these points of
intersection.

\bigskip
A second problem with real manifolds is that the number of
intersections may change with small perturbations of the map.  For
instance, if the map is changed so that $h(S)$ is tangent to $T$ and
has no other other intersections with $T$, then for small
perturbations $g$ near $h$, $g(S)$ may intersect $T$ in $0$, $1$, or
$2$ points.

\bigskip
Suppose now that $B$ is a bidisk in $\CC^2$, that $h$ is a complex
horseshoe map, and that $T$ and $S$ are complex submanifolds.  In this
case, there is a natural orientation on $T$ at any point given by
taking a vector $v$ in the tangent space of $T$ over this point and
using the set $\{v, i v\}$ to define the orientation at that point.
We can use the same 
procedure on $S$, then apply $h^n$ as before.  In this case, the
orientation induced on $\CC^2$ by $h^n(S)$ and $T$ is always the same
as the standard orientation.  Hence assigning $+1$ to such an
intersection and taking the sum gives the number of points in $T \cap
h^n(S)$.  

Additionally, if both $S$ and $T$ are complex manifolds, then the
number of intersections, counted with multiplicity, between $h^n(S)$
and $T$ is constant under small perturbations.  

Thus in studying question (i), we will use complex 1-dimensional
submanifolds.  

\bigskip
Recall that in the case of one variable, the Laplacian played a key
role by allowing us to relate the potential function $G$ to the 
measure $\mu$.  Here we make an extension of the Laplacian to $\CC^2$
in order to achieve a similar goal.  

For a function $f$ of two real variables $x$ and $y$, the exterior
derivative of $f$ is
$$ df = \frac{\partial f}{\partial x} dx + \frac{\partial f}{\partial
y} dy,  $$
which is invariant under smooth maps.  If we identify $\RR^2$ with
$\CC$ in the usual way, then multiplication by $i$ induces
the map $(i)^*$ on the cotangent bundle, and this map takes $dx$
to $dy$ and $dy$ to $-dx$.  Hence, defining $d^c = (i)^*d$, we have
$$ d^c f = \frac{\partial f}{\partial x} dy - \frac{\partial
f}{\partial y} dx, $$
which is invariant under smooth maps preserving the complex structure;
i.e., holomorphic maps.  Hence $d d^c$ is also invariant under
holomorphic maps. Expanding $d d^c$ gives
\begin{align*}
dd^c f & =  d\left(\frac{\partial f}{\partial x}\right) dy -
d \left(\frac{\partial f}{\partial y}\right) dx \\
& =  \left( \frac{\partial^2 f}{\partial x^2} + \frac{\partial^2
f}{\partial y^2} \right) dx dy, 
\end{align*}
which is nothing but the Laplacian.  This shows that the
Laplacian, when viewed as a map from functions to $2$-forms, is
invariant under holomorphic maps, and also shows that this procedure
can be carried out in any complex manifold of any dimension.
Moreover, it also shows that the property of being subharmonic is
invariant under holomorphic maps.

\begin{exercise}
Let $D_r$ be the disk of radius $r$ centered at $0$ in the plane, and
compute 
$$ \int_{D_r} d d^c \log|z|. $$
Hint:  This is equal to 
$$ \int_{\partial D_r} d^c \log|z|. $$
\end{exercise}

We next need to extend the idea of subharmonic functions to $\CC^2$.

\begin{definition}
A function $f: \CC^2 \ra \RR$ is {\bfseries plurisubharmonic} (psh) if $h$ is
u.s.c. and if the restriction of $h$ to any 1-dimensional complex line
satisfies the subaverage property.
\end{definition}

Intrinsically, an u.s.c. function $h$ is psh if and only if $dd^c h$
is nonnegative, where again we interpret this in the sense of
distributions.  

\bigskip
In fact, in the above definition we could replace the
phrase ``complex line'' by ``complex submanifold'' without changing
the class of functions, since subharmonic functions are invariant
under holomorphic maps.  As an example of the usefulness of this and
the invariance property of $dd^c$, suppose that $\phi$ is a
holomorphic embedding of $\CC$ into $\CC^2$ and that $h$ is smooth and
psh on $\CC^2$.  Then we can either evaluate $dd^c h$ and pull back
using $\phi$, or we can first pull back and then apply $dd^c$.  In
both cases we get the same positive distribution on $\CC$.

\section{Currents and applications to dynamics}

In this section we provide a brief introduction to the theory of
currents.  A current is simply a linear functional on the space of
smooth differential forms, and hence may be viewed as a generalization
of measures.  That is, a current $\mu$ acts on a differential form of
a given degree, say $\phi = f_1 dx + f_2 dy$ in the case of a $1$-form,
to give a complex number $\mu(\phi)$, and this assignment is linear in
$\phi$.  This is a generalization of a measure in the sense that a
measure acts on $0$-forms (functions) by integrating the function
against the measure.  

As an example, suppose that $M \subseteq \CC^2$ is a submanifold of
real dimension $n$.  Then the current of integration associated to $M$
is a current acting on $n$-forms $\phi$, and is simply given by
$$ [M](\phi) = \int_M \phi.  $$
In this example the linearity is immediate, as is the relationship to
measures.  Note that in particular, if $p \in \CC^2$, then $[p] =
\delta_p$, the delta mass at $p$, acts on $0$-forms.

\begin{example}
Suppose $P : \CC \ra \CC$ is a polynomial having only simple roots,
and let $R$ be the set of roots of $P$.  Then $[R]$ is a current
acting on $0$-forms, and 
$$ [R] = \frac{1}{2 \pi} d d^c \log|P|.  $$
This formula is still true for arbitrary polynomials if we account for
multiplicities in constructing $[R]$.
\end{example}

We can extend this last example to the case of polynomials from
$\CC^2$ to $\CC$.  This is the content of the next proposition.

\begin{prop}
If $P : \CC^2 \ra \CC$ is a polynomial and $V = \{ P = 0 \}$, then 
$$ [V] = \frac{1}{2\pi} d d^c \log|P|, $$
where again $[V]$ is interpreted with weights according to
multiplicity.  (This is known as the Poincar\'e-Lelong formula.)
\end{prop}

Suppose now that $f : \CC^2 \ra \CC^2$ is a H\'enon diffeomorphism,
and define
\begin{align*}
G^+(p) & = \lim_{n \ra \infty} \frac{1}{2^n} \log^+|\pi_1(f^n(p))|
\\
G^-(p) & = \lim_{n \ra \infty} \frac{1}{2^n}
\log^+|\pi_2(f^{-n}(p))|,
\end{align*}
where $\pi_j$ is projection to the $j$th coordinate.  As in the case
of the function $G$ defined for a 1-variable polynomial, it is not
hard to check that $G^+$ is psh, is identically $0$ on $K^+$, is
pluriharmonic on $\CC^2 - K^+$ (i.e., is harmonic on any complex
line), and is positive on $\CC^2 - K^+$.  In analogy with the function
$G$, we say that $G^+$ is the Green function of $K^+$.  Likewise,
$G^-$ is the Green function of $K^-$.  

Note that
for $n$ large and $p \not \in K$, $|\pi_1 f^n(p)|$ is comparable to
the square of $|\pi_2 f^n(p)|$, and hence we may replace $|\pi_1
f^n(p)|$ by $\|f^n(p)\|$ in the formula for $G^+$, and likewise for
$G^-$.  

\bigskip
Again in analogy with the 1-variable case, and using the equivalence
between the Laplacian and $d d^c$ outlined earlier, we define
\begin{align*}
\mu^+ & = \frac{1}{2 \pi} d d^c G^+ \\
\mu^- & = \frac{1}{2 \pi} d d^c G^-.
\end{align*}
Then $\mu^+$ and $\mu^-$ are currents supported on $J^+ = \partial
K^+$ and $J^- = \partial K^-$, respectively.  Moreover, $\mu^\pm$
restrict to measures on complex 1-dimensional submanifolds in the
sense that we can pull back $G^\pm$ from the submanifold to an open
set in the plane, then take $d d^c$ on this open set.  

\bigskip
As an analog of part 1 of Brolin's theorem, we have the following
theorem.

\begin{theorem}
Let $V$ be the (complex) $x$-axis in $\CC^2$, i.e., the set where
$\pi_2$ vanishes, and let $f$ be a complex H\'enon map.  Then
$$ \lim_{n \ra \infty} \frac{1}{2^n} [f^{-n} V] =
\mu^+.  $$
\end{theorem}

{\bf Proof:}  Note that the set $f^{-n}V$ is the set where the
polynomial $\pi_1 f^n$ vanishes.  Hence the previous proposition
implies that
$$ [f^{-n} V] = \frac{1}{2\pi} d d^c \log|\pi_1 f^n|.  $$
Passing to the limit and using an argument like that in Brolin's
theorem to replace $\log$ by $\log^+$, we obtain the theorem.  See
\cite{bs1} or \cite{fs1} for more details.  $\BOX$

\bigskip
As a more comprehensive form of this theorem, we have the following.

\begin{theorem}  \label{thm:converge}
If $S$ is a complex disk and $f$ is a complex H\'enon map, then
$$ \lim_{n \ra \infty} \frac{1}{2^n} [f^{-n} S] = c \mu^+, $$
where $c = \mu^-[S]$.  An analogous statement is true with $\mu^+$ and
$\mu^-$ interchanged and $f^n$ in place of $f^{-n}$.
\end{theorem}

As a corollary, we obtain the following theorem from \cite{bs3}.

\begin{cor}
If $p$ is a periodic saddle point, then $W^u(p)$ is dense in $J^-$.
\end{cor}

{\bf Proof:} Replacing $f$ by $f^n$, we may assume that $p$ is fixed.
Let $S$ be a small disk in $W^u(p)$ containing $p$.  Then 
the forward iterates of $S$ fill out the entire unstable manifold.
Moreover, by the previous theorem, the currents associated with these
iterates converge to $c \mu^-$ where $c = \mu^+[S]$.  If $c \neq 0$,
then the corollary is complete since then $W^u(p)$ must be dense in 
${\rm supp} \; \mu^- = J^-$.  

But $c$ cannot be $0$, because if it were, then $G^+|S$ would be
harmonic, hence identically $0$ by the minimum principle since $G$ is
nonnegative on $S$ and $0$ at $p$.  Hence $S$ would be contained in
$K$, which is impossible since the iterates of $S$ fill out all of
$W^u(p)$, which is not bounded.   Thus $c \neq 0$, so $W^u(p)$ is
dense in $J^-$.  $\BOX$

\bigskip
This corollary gives some indication of why pictures of invariant sets
on complex slices in $\CC^2$ show essentially the full complexity of
the map.  If we start with any complex slice which is transverse to
the stable manifold of a periodic point $p$, then the forward iterates of
this slice accumulate on the unstable manifold of $p$, hence on all of
$J^-$ by the corollary.  All of this structure is then reflected in
the original slice, giving rise to sets which are often self-similar
and bear a striking resemblance to Julia sets in the plane.  

\section{Currents and H\'enon maps}

In this section we continue the study of the currents $\mu^+$ and
$\mu^-$ in order to obtain more dynamical information.

We first consider this in the context of the horseshoe map.  Recall
that $B$ is a square in the plane and that we have defined measures
$m^+$ and $m^-$, and their product measure $m$ in
section~\ref{sec:m-}.  

In fact, $m^+$ and $m^-$ generalize to $\mu^+$
and $\mu^-$ in the case that the H\'enon map is a horseshoe.  More
explicitly, let 
$D_\lambda$ be a family of complex disks in $\CC^2$ which intersect
$\RR^2$ in a horizontal segment in $B$ and such that these segments
fill out all of $B$.  Then we can recover $\mu^+$, at least restricted
to $B$, by
$$ \mu^+|B = \int [D_\lambda] d m^+(\lambda).  $$

\bigskip
In analogy with the construction of $m$ as a product measure using
$m^+$ and $m^-$, we would like to combine $\mu^+$ and $\mu^-$ to
obtain a measure $\mu$.  Since $\mu^+$ and $\mu^-$ are currents, the
natural procedure to try is to take $\mu = \mu^+ \wedge \mu^-$.  While
forming the wedge product is not well-defined for arbitrary currents,
it is well-defined in this case using the fact that these currents are
obtained by taking $d d^c$ of a continuous psh function and applying a
theorem of pluripotential theory.  In this way we get a measure $\mu$
on $\CC^2$. 

\begin{definition} 
$\mu = \mu^+ \wedge \mu^-$.
\end{definition}

\bigskip
We next collect some useful facts about $\mu$. 

\paragraph{(1)} $\mu$ is a probability measure.  For a proof of
this, see \cite{bs1}.

\paragraph{(2)} $\mu$ is invariant under $f$.  To see this, note that
since 
$$ G^{\pm} = \lim_{n \ra \infty} \frac{1}{2^n} \log^+ \| f^{\pm n}\|, $$
we have $G^\pm(f(p)) = 2^\pm G^\pm(p)$.  Since $\mu^\pm = (1/2\pi) d
d^c G^\pm$, this implies that $f^*(\mu^\pm) = 2^\pm \mu^\pm$, and
hence
\begin{align*}
f^*(\mu) & = f^*(\mu^+) \wedge f^*(\mu^-) \\
& = 2 \mu^+ \wedge \frac{1}{2} \mu^- \\
& = \mu^+ \wedge \mu^-\\
& = \mu. 
\end{align*}

\begin{definition}
$ J = J^+ \cap J^-$.
\end{definition}

\paragraph{(3)} ${\rm supp}(\mu) \subseteq J$.  This is a simple
consequence of the fact that the support of $\mu$ is contained in the
intersection of ${\rm supp}(\mu^+) = J^+$ and ${\rm supp}(\mu^-) =
J^-$ and the definition of $J$.

\bigskip
In order to examine the support of $\mu$ more precisely, we turn our
attention for a moment to Shilov boundaries.  Let $X$ be a subset
either of $\CC$ or $\CC^2$.  We say that a set $B$ is a {\em boundary} for
$X$ if $B$ is closed and if for any holomorphic polynomial
$P$ we have
$$ \max_X |P| = \max_B |P|. $$
With the right conditions, the intersections of any set of boundaries
is a again a boundary by a theorem of Shilov, so we can intersect them
all to obtain the smallest such boundary.  This is called the {\em\bfseries
Shilov boundary} for $X$.  

\begin{example}
Let $X = D_1 \times D_1$, where $D_1$ is the unit disk.  Then the
Shilov boundary for $X$ is $(\partial D_1) \times (\partial D_1)$,
while the topological boundary for $X$ is 
$$ \partial X = (D_1 \times \partial D_1) \cup (\partial D_1 \times
D_1). $$
\end{example}

The following theorem is contained in \cite{bt}.  
 
\begin{theorem}
${\rm supp}(\mu) = \partial_{\rm Shilov}(K)$. 
\end{theorem}

We have already defined $J$ as the intersection of $J^+$ and $J^-$,
and the choice of notation is designed to suggest an analogy with the
Julia set in one variable.  However, in 2 variables, the support of
$\mu$ is also a natural candidate as a kind of Julia set.  Hence we
make the following definition.

\begin{definition}
$J^* = {\rm supp}(\mu)$.
\end{definition}

\section{Heteroclinic points and Pesin theory}

In the previous section, we discussed some of the formal properties of
$\mu$ arising from considerations of the definition and of potential
theory.  In this section we concentrate on the less formal properties
of $\mu$ and on the relation of $\mu$ to periodic points.  The philosophy here
is that since $\mu^+$ and $\mu^-$ describe the distribution of
$1$-dimensional objects, $\mu$ should describe the distribution of
$0$-dimensional objects.

\bigskip
An example of a question using this philosophy is the following.  For
a periodic point $p$, we know that $\mu^+$ describes the distribution
of $W^s(p)$ and $\mu^-$ describes the distribution of $W^u(p)$.  Does
$\mu$ describe (in some sense) the distribution of intersections
$W^s(p) \cap W^u(q)$?

\begin{definition}  Let $p$ and $q$ be saddle periodic points of a
diffeomorphism $f$.  A point in the set $(W^u(p) \cap W^s(q)) -
\{p,q\}$ is called a {\bfseries heteroclinic point}.  If $p=q$, then
such a point is called a {\bfseries homoclinic point}.  
\end{definition}

Unfortunately, none of the techniques discussed so far allow us to
show that there is even one heteroclinic point.  That is, it is
possible for the stable manifold of $p$ and the unstable manifold of
$q$ to have an empty intersection.  

To see how this could happen, we
first need the fact that for a holomorphic diffeomorphism $f$ of
$\CC^2$ with a saddle point $p$, there is an injective holomorphic map
$\phi_u: \CC \ra W^u(p)$ which maps {\em onto} the unstable manifold,
and likewise for the stable manifold.  

Now, if $\pi_j$ represents projection onto the $j$th coordinate, then
$\pi_1 \phi_u : \CC \ra \CC$ is an entire function, and as such can
have an omitted value.  As an example, $\pi_1 \phi_u(z)$ could be equal
to $e^z$ and hence would omit the value $0$.  It could happen that
there is a second saddle point $q$ such that $W^s(q)$ is the $x$-axis,
in which case $W^u(p) \cap W^s(q) = \emptyset$.  

\bigskip
At first glance, it may seem that this contradicts some of our earlier
results.  It might seem that theorem~\ref{thm:converge} should imply
that $W^u(p)$ intersects transversals which cross $J^-$, but in fact,
that statement is a statement about convergence of distributions.
Each of the distributions must be evaluated against a test function,
and the test function must be positive on an open set.  Thus there is
still room for $W^u(p)$ and $W^s(p)$ to be disjoint.

\bigskip
Hence, in order to understand more about heteroclinic points, we need
a better understanding of the stable and unstable manifolds.  One
possible approach is to use what is known as Ahlfors' three island
theorem.  This theorem concerns entire maps $\psi:\CC \ra \CC$.
Roughly, it says that if we have $n$ open regions in the plane and
consider their inverse images under $\psi$, then some fixed proportion
of them will have a preimage which is compact and which maps
injectively under $\psi$ onto the corresponding original region.  

If we apply this theorem to the map $\pi_1 \phi_u$ giving $W^u(p)$,
then we can divide the plane into increasingly more and smaller
islands, and we can do this in such a way that at each stage we obtain
more of $W^u(p)$ as the injective image of regions in the plane.  The
result is that we get a picture of $W^u(p)$ which is locally laminar.

Since $W^u(p)$ is dense in $J^-$, this gives us one possible approach
to studying $\mu^-$, and we can use a similar procedure to study
$\mu^+$.  However, recall that our goal here is to describe
heteroclinic points.  Thus in order for this approach to apply, we
need to be able to get the disks for $\mu^+$ to intersect the disks
for $\mu^-$.  Unfortunately, we don't get any kind of uniformity in
the disks using this approach, so getting this intersection is
difficult.  

\bigskip
An alternative approach is to use the theory of hyperbolicity in the
sense of Oseledec and Pesin.  Since $\mu$ is a (non-uniformly)
hyperbolic measure with respect to this theory, we get that at $\mu$-almost
every point of $\CC^2$ there are stable and unstable manifolds and that
these manifolds are transverse.  We can then identify the stable and
unstable manifolds obtained using this theory with the disks obtained
in the previous non-uniform laminar picture to guarantee that we get
intersections between stable and unstable manifolds and hence
heteroclinic points.  Putting all of this together, we obtain the
following theorem, contained in \cite{bls1}.

\begin{theorem}   \label{thm:jstar}
$J^*$ is the closure of
\begin{enumerate}
\item the set of all periodic saddle points,
\item the union of all $W^u(p) \cap W^s(q)$ over all periodic saddles
$p$ and $q$.
\end{enumerate}
\end{theorem}

This theorem can be viewed as an analog of the theorem in one variable
dynamics that the Julia set is the closure of the repelling periodic
points.  For this reason, the set $J^*$ is perhaps a better analogue
of the Julia set in the two dimensional case than is $J$.

Recall that $J^* = \partial_{\rm shilov} K \subseteq \partial K = J$.
In the case that $f$ is an Axiom A diffeomorphism, it is a theorem
that $J^* = J$.  However, it is an interesting open question whether
this equality 
holds in general.  If it were the case that $J \neq J^*$, then
there would be a saddle periodic point $q$ and another point $p$ such
that $p \in \overline{W^s(q)} \cap \overline{W^u(q)}$, but $p \not \in
\overline{W^s(q) \cap W^u(q)}$. 

\bigskip
In fact, using the ideas of Pesin theory, one can get
precise information about the number of periodic points of a given
period and how their distribution relates to the measure $\mu$.
This is contained in the following theorem and corollary, contained in
\cite{bls}. 

\begin{theorem}
Let $f: \CC^2 \ra \CC^2$ be a complex H\'enon map, and let $P_n$ be
either the set of fixed points of $f^n$ or the set of saddle points of
minimal period $n$.  Then
$$ \lim_{n \ra \infty} \frac{1}{2^n} \sum_{p \in P_n} \delta_p =
\mu.$$
\end{theorem}

For the following corollary, let $P_n$ be the set of saddle points of
$f$ of minimal period $n$, and let $|P_n|$ denote the number of points
contained in this set.  

\begin{cor}
There are periodic saddle points of all but finitely many periods.  More
precisely, we have
$$ \lim_{n \ra \infty} \frac{|P_n|}{2^n} = 1. $$
\end{cor}

Recall that the horseshoe map had periodic points of all
periods, so while we haven't achieved that result for general H\'enon
maps, we have still obtained a good deal of information about periodic
points and heteroclinic points.

\section{Topological entropy}

Recall that the horseshoe map is topologically equivalent to the shift
map on two symbols.  One could also ask if it is topologically
equivalent to the shift on four symbols.  That is, if $h$ is the
horseshoe map defined on the square $B$, then $h(B) \cap h^{-1}(B)$
consists of four components, and we can label these components with
four symbols.  However, with this labeling scheme, one can check by
counting that
not all sequences of symbols correspond to an orbit of a point in the
way that sequences of two symbols did.  In fact, the number of symbol
sequences of length 2 corresponding to part of an orbit is 8, while
the total number of possible sequences of length 3 is 16.  Allowing
longer sequences and letting $S(n)$ denote the number of sequences of
length $n$ which correspond to part of an orbit, we obtain the formula
$$ \lim_{n \ra \infty} \frac{1}{n} \log S(n) = \log 2. $$
The number $\log 2$ is the topological entropy of the horseshoe map,
and can be defined as the maximum growth rate over all finite
partitions.  In general, the shift map on $N$ symbols has entropy
$\log N$, and since entropy is a topological invariant, we see that
all of these different shift maps are topologically distinct.

In the case of a general H\'enon map, we have the following theorem,
contained in \cite{smillie}.

\begin{theorem}
The topological entropy of a complex H\'enon map is $\log 2$.
\end{theorem}
 
Topological entropy is a useful idea because it is connected to many
different aspects of polynomial automorphisms.  It is a measure of area
growth and of the growth rate of the number of periodic points, both
of which are closely related to the degree of the map as a
polynomial.  Moreover, it is related to measure theoretic entropy in
the sense that for any probability measure $\nu$, the measure
theoretic entropy, $h_\nu(f)$, is bounded from above by the
topological entropy $h_{\rm top}(f)$.  Moreover, $\mu$ is the {\em
unique} measure for which $h_\mu(f) = h_{\rm top}(f)$.  

\bigskip
We can also consider topological entropy for real H\'enon maps.  In
contrast to the theorem above, in this case we have $0 \leq h_{\rm
top}(f_{\bf R}) \leq \log 2$, and all values are possible.  However,
one can show 
that not all values are possible for Axiom A automorphisms, but only
logarithms of algebraic numbers \cite{f}.  Moreover, we also have the
following theorem  \cite{bls1}.

\begin{theorem}   \label{thm:conds}
For a real H\'enon map $f_{\bf R}$, the following are equivalent.

\begin{enumerate}
\item $h_{\rm top}(f_{\bf R}) = \log 2$.
\item $J^* \subseteq \RR^2$.
\item $K \subseteq \RR^2$.
\item All periodic points are real.
\end{enumerate}

Moreover, these conditions imply that $J = J^*$.
\end{theorem}

{\bf Proof:}  Condition 1 implies that $f_{\bf R}$ has a measure
$\mu'$ of maximal entropy with $\text{supp}(\mu') \subseteq \RR^2$.
By uniqueness we have $\mu' = \mu^*$, so $\text{supp}(\mu^*) \subseteq
\RR^2$, thus giving condition 2.

Condition 2 implies that $J^* = \partial_{\rm Shilov} K$ is contained in
$\RR^2$, which implies that $K$ is contained in $\RR^2$.  This gives
condition 3, and in fact, since polynomials in $\RR^2$ are dense in
the set of continuous functions of $\RR^2$, this also implies that
$\partial_{\rm Shilov} K = K$, and hence $J^*=K$ and thus $J^*=J$
since $J^* \subseteq J \subseteq K$.

Condition 3 immediately implies condition 4.

Condition 4 together with theorem~\ref{thm:jstar} implies that $J^*
\subseteq \RR^2$, which implies that $\text{supp}(\mu^*) \subseteq
\RR^2$, which implies condition 1. 
$\BOX$ 

\bigskip
These conditions are true for the set of
real H\'enon maps which are horseshoes.  We can identify such maps
with their parameter values in $\RR^2$, in which case the set of
horseshoe maps is an open set in $\RR^2$.  Since topological entropy
is continuous for $C^\infty$ diffeomorphisms, we see that maps on the
boundary of this set also satisfy the above conditions.

\bigskip
These conditions also apply in the following theorem, from \cite{bs5}.
Recall that a homoclinic intersection is an intersection of $W^u(p)$
and $W^s(p)$ at some point $q \neq p$ for some saddle point $p$.  This
intersection is a homoclinic tangency if the stable and unstable
manifolds are tangent at $q$, and this is a quadratic tangency if the
manifolds have quadratic contact at $q$.

\begin{theorem}
If the previous conditions hold, then 
\begin{enumerate}
\item periodic points are dense in $K$,
\item every periodic point is a saddle with expansion constants
bounded below,
\item either $f$ is Axiom A or $f$ has a quadratic 
homoclinic tangency.
\end{enumerate}
\end{theorem}

Suppose we have a 1-parameter family of real H\'enon maps which starts out as a
horseshoe, then passes through a homoclinic tangency, so that as we
increase the parameter value, some local pieces of the stable and unstable
manifold in $\RR^2$ first intersect in 2 points, then at one tangent
point, then 
don't intersect at all.  In this case, the
intersections of the stable and unstable manifolds move out of
$\RR^2$, which causes a decrease in topological entropy by
theorem~\ref{thm:conds} and the fact that $\log 2$ is the maximum
possible entropy for a real H\'enon map.  Since topological entropy is
continuous and is a topological invariant, it follows that we pass
through infinitely many conjugacy classes as we change the parameter.
That is, for two maps with different topological entropy, there is no
homeomorphism which conjugates one to the other.  This presents a
striking contrast to the horseshoe example in which small changes of
the original horseshoe were all conjugate to the shift map on two symbols.

\section{Conclusion}

We have presented here an overview of some of the major techniques and
results in the study of the iteration of polynomial automorphisms of
$\CC^2$ and have demonstrated some of the influences coming from real
and measurable dynamics as well as complex dynamics of one variable.  

There has also been a great deal of work in other directions in the
study of complex dynamics in several variables, including the study of
rational maps on complex projective space by Forn{\ae}ss, Sibony, Ueda,
and others, and the
study of holomorphic vector fields by Ahern, Bass, Coomes, Forn{\ae}ss,
Forstneric, Grellier, Meisters, Sibony, Suzuki, and others.  There is
also work on nonpolynomial automorphisms of $\CC^2$ by Buzzard,
Forn{\ae}ss, Sibony, and others, as well as the study of foliations of
$\PP^2$ by Cano, Camacho, Gomez-Mont, Lins-Neto, Sad, and others.
There are many open problems in each of these areas.

\footnotesize
\sc
\noindent
John Smillie,
Department of Mathematics,
Cornell University,
Ithaca, NY  14853 
({\tt smillie@@math.cornell.edu})

\medskip\noindent
Gregery T. Buzzard,
Department of Mathematics,
Indiana University,
Bloomington, IN  47405
({\tt gbuzzard@@msri.org, gbuzzard@@indiana.edu})

\end{document}